\documentclass[12pt]{article}
\usepackage{amssymb}

\newtheorem{definition}{Definition}[section]
\newtheorem{theorem}[definition]{Theorem}

\newtheorem{lemma}[definition]{Lemma}

\newtheorem{remark}[definition]{Remark}
\newtheorem{corollary}[definition]{Corollary}
\newtheorem{problem}[definition]{Problem}
\def\R{\mathbb R}
\def\C{\mathbb C}

\newcommand{\beast}{\begin{eqnarray*}}
\newcommand{\eeast}{\end{eqnarray*}}

\begin{document}

\title{
The displacement and split decompositions \\
for a 
$Q$-polynomial 
distance-regular
graph
\author{Paul Terwilliger}
}
\date{}
\maketitle

\begin{abstract}
Let $\Gamma$ denote a $Q$-polynomial
distance-regular graph
with  diameter at least three 
and standard module $V$.
We introduce two direct sum decompositions of $V$.
We call these the {\it displacement} decomposition for $\Gamma$
and the {\it split} decomposition for $\Gamma$.
We describe how
these decompositions are related.
\end{abstract}

\section{Introduction}
In this paper $\Gamma=(X,R)$ will denote a $Q$-polynomial
distance-regular graph
with  diameter 
$D \geq 3$ and adjacency matrix $A$ (see Section 2 for formal definitions).
In order to describe our main results we make a few comments.
Fix a vertex $x \in X$. For $0 \leq i \leq D$ let
 $E^*_i=E^*_i(x) 
$ denote the diagonal matrix in $\mbox{Mat}_X(\C)$
that represents the
 projection onto the $i$th subconstituent of
$\Gamma$ with respect to $x$.
Let $E_0, E_1, \ldots, E_D$
denote a $Q$-polynomial ordering of the primitive idempotents for $A$
and let
$A^*=A^*(x)$ denote the corresponding dual adjacency 
matrix. 
The {\it subconstituent algebra}  $T=T(x)$ is the subalgebra of
$\mbox{Mat}_X(\C)$
generated by $A$ and $A^*$.
Let $W$ denote an irreducible $T$-module.
By the {\it displacement} of $W$ we mean
$\rho+\tau+d-D$, where
$\rho=\hbox{min}\lbrace i |E^*_iW\not=0 \rbrace $,
$\tau=\hbox{min}\lbrace i |E_iW\not=0 \rbrace $,
$d=|\lbrace i |E_iW\not=0 \rbrace|-1 $.
We show the displacement of $W$ is nonnegative
and at most $D$.
Let $V=\C^X$ denote the standard module.
We show
$V=\sum_{\eta=0}^DV_\eta$ (orthogonal direct sum),
where
$V_\eta$ denotes the
subspace of $V$
spanned by the irreducible $T$-modules
that have displacement $\eta$. 
This is the {\it displacement decomposition} with respect to $x$.
For $-1\leq i,j\leq D$ we define
$V_{ij} = (E^*_0V+\cdots+E^*_iV)\cap (E_0V+\cdots +E_jV)
$.
We show
$V = \sum_{i=0}^D\sum_{j=0}^D {\tilde V}_{ij}$ (direct sum),
where ${\tilde V}_{ij}$ denotes 
the orthogonal complement of $V_{i,j-1}+V_{i-1,j}$ in
$V_{ij}$ with respect to the Hermitean dot product.
This direct sum is the {\it split decomposition} with respect to $x$.
The above
decompositions are related as follows.
For $0\leq \eta \leq D$ we show
$V_\eta= 
\sum {\tilde V}_{ij}$,
where the sum is over all ordered pairs $i,j$
such that $0 \leq i,j\leq D$ and $i+j=D+\eta$.
Using this we obtain the following results.
For $0 \leq i,j\leq D$ we show 
$V_{ij}=0$ if $i+j<D$.
For $0 \leq i\leq D$ let $\theta_i$ (resp. $\theta^*_i$)
denote the eigenvalue
of $A$ (resp. $A^*$) for $E_i$ (resp. $E^*_i$).
For $0 \leq i,j\leq D$ we show
 $(A-\theta_jI){\tilde V}_{ij}
\subseteq {\tilde V}_{i+1,j-1}$
and
$(A^*-\theta^*_iI){\tilde V}_{ij}
\subseteq {\tilde V}_{i-1,j+1}$,
where
${\tilde V}_{rs}:=0$ unless $r,s \in \lbrace 0,1,\ldots, D\rbrace $.
We finish with an application related to the work of Brouwer,
Godsil, Koolen and Martin
\cite{wdw}
concerning the dual width of a subset of $X$.

\section{Preliminaries concerning distance-regular graphs}
In this section 
we review some definitions and basic concepts concerning distance-regular
graphs.
For more background information we refer the reader to 
\cite{bannai}, \cite{bcn}, \cite{godsil} and \cite{terwSub1}.

\medskip
\noindent
Let $\C$ denote the complex number field.
Let $X$ denote a nonempty  finite  set.
Let
 $\hbox{Mat}_X(\C)$ 
denote the $\C$-algebra
consisting of all matrices whose rows and columns are indexed by $X$
and whose entries are in $\C  $. Let
$V=\C^X$ denote the vector space over $\C$
consisting of column vectors whose 
coordinates are indexed by $X$ and whose entries are
in $\C$.
We observe
$\hbox{Mat}_X(\C)$ 
acts on $V$ by left multiplication.
We call $V$ the {\it standard module}.
We endow $V$ with the Hermitean inner product $\langle \, , \, \rangle$ 
that satisfies
$\langle u,v \rangle = u^t\overline{v}$ for 
$u,v \in V$,
where $t$ denotes transpose and $\overline{\phantom{v}}$
denotes complex conjugation.
For all $y \in X,$ let $\hat{y}$ denote the element
of $V$ with a 1 in the $y$ coordinate and 0 in all other coordinates.
We observe $\{\hat{y}\;|\;y \in X\}$ is an orthonormal basis for $V.$

\medskip
\noindent
Let $\Gamma = (X,R)$ denote a finite, undirected, connected graph,
without loops or multiple edges, with vertex set $X$ and 
edge set
$R$.   
Let $\partial $ denote the
path-length distance function for $\Gamma $,  and set
$D := \mbox{max}\{\partial(x,y) \;|\; x,y \in X\}$.  
We call $D$  the {\it diameter} of $\Gamma $.
 We say $\Gamma$ is {\it distance-regular}
whenever for all integers $h,i,j\;(0 \le h,i,j \le D)$ 
and for all
vertices $x,y \in X$ with $\partial(x,y)=h,$ the number
\begin{eqnarray*}
p_{ij}^h = |\{z \in X \; |\; \partial(x,z)=i, \partial(z,y)=j \}|
\end{eqnarray*}
is independent of $x$ and $y.$ The $p_{ij}^h$ are called
the {\it intersection numbers} of $\Gamma.$ 

\medskip
\noindent
For the rest of this paper we assume  $\Gamma$  
is  distance-regular  with diameter $D\geq 3$. 

\medskip
\noindent 
We recall the Bose-Mesner algebra of $\Gamma.$ 
For 
$0 \le i \le D$ let $A_i$ denote the matrix in $\hbox{Mat}_X(\C)$ with
$xy$ entry
$$
{(A_i)_{xy} = \cases{1, & if $\partial(x,y)=i$\cr
0, & if $\partial(x,y) \ne i$\cr}} \qquad (x,y \in X).
$$
We call $A_i$ the $i$th {\it distance matrix} of $\Gamma.$
We abbreviate $A:=A_1$ and call this the {\it adjacency
matrix} of $\Gamma.$ We observe
(i) $A_0 = I$;
 (ii)
$\sum_{i=0}^D A_i = J$;
(iii)
$\overline{A_i} = A_i \;(0 \le i \le D)$;
(iv) $A_i^t = A_i  \;(0 \le i \le D)$;
(v) $A_iA_j = \sum_{h=0}^D p_{ij}^h A_h \;( 0 \le i,j \le D)
$,
where $I$ (resp. $J$) denotes the identity matrix 
(resp. all 1's matrix) in 
 $\hbox{Mat}_X(\C)$.
 Using these facts  we find
 $A_0,A_1,\ldots,A_D$
is a basis for a commutative subalgebra $M$ of 
$\mbox{Mat}_X(\C)$.
We call $M$ the {\it Bose-Mesner algebra} of $\Gamma$.
It turns out $A$ generates $M$ \cite[p.~190]{bannai}.
By \cite[p.~45]{bcn}, $M$ has a second basis 
$E_0,E_1,\ldots,E_D$ such that
(i) $E_0 = |X|^{-1}J$;
(ii) $\sum_{i=0}^D E_i = I$;
(iii) $\overline{E_i} = E_i \;(0 \le i \le D)$;
(iv) $E_i^t =E_i  \;(0 \le i \le D)$;
(v) $E_iE_j =\delta_{ij}E_i  \;(0 \le i,j \le D)$.
We call $E_0, E_1, \ldots, E_D $  the {\it primitive idempotents}
of $\Gamma$.  

\medskip
\noindent
We  recall the eigenvalues
of  $\Gamma $.
Since $E_0,E_1,\ldots,E_D$ form a basis for  
$M$ there exist complex scalars $\theta_0,\theta_1,
\ldots,\theta_D$ such that
$A = \sum_{i=0}^D \theta_iE_i$.
Observe
$AE_i = E_iA =  \theta_iE_i$ for $0 \leq i \leq D$.
By \cite[p.~197]{bannai} the 
scalars $\theta_0,\theta_1,\ldots,\theta_D$ are
in $\R.$ Observe
$\theta_0,\theta_1,\ldots,\theta_D$ are mutually distinct 
since $A$ generates $M$. We call $\theta_i$  the {\it eigenvalue}
of $\Gamma$ associated with $E_i$ $(0 \leq i \leq D)$.
Observe 
\begin{eqnarray*}
V = E_0V+E_1V+ \cdots +E_DV \qquad \qquad {\rm (orthogonal\ direct\ sum}).
\end{eqnarray*}
For $0 \le i \le D$ the space $E_iV$ is the  eigenspace of $A$ associated 
with $\theta_i$.

\medskip
\noindent 
We now recall the Krein parameters.
Let $\circ $ denote the entrywise product in
$\mbox{Mat}_X(\C)$.
Observe
$A_i\circ A_j= \delta_{ij}A_i$ for $0 \leq i,j\leq D$,
so
$M$ is closed under
$\circ$. Thus there exist complex scalars
$q^h_{ij}$  $(0 \leq h,i,j\leq D)$ such
that
$$
E_i\circ E_j = |X|^{-1}\sum_{h=0}^D q^h_{ij}E_h
\qquad (0 \leq i,j\leq D).
$$
By \cite[p.~170]{Biggs}, 
$q^h_{ij}$ is real and nonnegative  for $0 \leq h,i,j\leq D$.
The $q^h_{ij}$ are called the {\it Krein parameters}.
The graph $\Gamma$ is said to be {\it $Q$-polynomial}
(with respect to the given ordering $E_0, E_1, \ldots, E_D$
of the primitive idempotents)
whenever for $0 \leq h,i,j\leq D$, 
$q^h_{ij}= 0$
(resp. 
$q^h_{ij}\not= 0$) whenever one of $h,i,j$ is greater than
(resp. equal to) the sum of the other two
\cite{bannai,
wdw,
caugh1,
caugh2,
curtin3,
curtin4,
dickie1,
dickie2,
aap1,
aap2}.
From now on assume $\Gamma$ is 
$Q$-polynomial with respect to $E_0, E_1, \ldots, E_D$.

\medskip
\noindent
We  recall the dual Bose-Mesner algebra of $\Gamma.$
Fix a vertex $x \in X.$ We view $x$ as a ``base vertex.''
For 
$ 0 \le i \le D$ let $E_i^*=E_i^*(x)$ denote the diagonal
matrix in $\hbox{Mat}_X(\C)$ with $yy$ entry
\begin{equation}\label{DEFDEI}
{(E_i^*)_{yy} = \cases{1, & if $\partial(x,y)=i$\cr
0, & if $\partial(x,y) \ne i$\cr}} \qquad (y \in X).
\end{equation}
We call $E_i^*$ the  $i$th {\it dual idempotent} of $\Gamma$
 with respect to $x$ \cite[p.~378]{terwSub1}.
We observe
(i) $\sum_{i=0}^D E_i^*=I$;
(ii) $\overline{E_i^*} = E_i^*$ $(0 \le i \le D)$;
(iii) $E_i^{*t} = E_i^*$ $(0 \le i \le D)$;
(iv) $E_i^*E_j^* = \delta_{ij}E_i^* $ $(0 \le i,j \le D)$.
By these facts 
$E_0^*,E_1^*, \ldots, E_D^*$ form a 
basis for a commutative subalgebra
$M^*=M^*(x)$ of 
$\hbox{Mat}_X(\C).$ 
We call 
$M^*$ the {\it dual Bose-Mesner algebra} of
$\Gamma$ with respect to $x$ \cite[p.~378]{terwSub1}.
For $0 \leq i \leq D$ let $A^*_i = A^*_i(x)$ denote the diagonal
matrix in 
 $\hbox{Mat}_X(\C)$
with $yy$ entry
$(A^*_i)_{yy}=\vert X \vert (E_i)_{xy}$ for $y \in X$.
Then $A^*_0, A^*_1, \ldots, A^*_D$ is a basis for $M^*$ 
\cite[p.~379]{terwSub1}.
Moreover
(i) $A^*_0 = I$;
(ii)
$\overline{A^*_i} = A^*_i \;(0 \le i \le D)$;
(iii) $A^{*t}_i = A^*_i  \;(0 \le i \le D)$;
(iv) $A^*_iA^*_j = \sum_{h=0}^D q_{ij}^h A^*_h \;( 0 \le i,j \le D)
$
\cite[p.~379]{terwSub1}.
We call 
 $A^*_0, A^*_1, \ldots, A^*_D$
the {\it dual distance matrices} of $\Gamma$ with respect to $x$.
We abbreviate
$A^*:=A^*_1$ 
and call this the {\it dual adjacency matrix} of $\Gamma$ with
respect to $x$.
The matrix $A^*$ generates $M^*$ \cite[Lemma 3.11]{terwSub1}.

\medskip
\noindent We recall the dual eigenvalues of $\Gamma$.
Since $E^*_0,E^*_1,\ldots,E^*_D$ form a basis for  
$M^*$, there exist complex scalars $\theta^*_0,\theta^*_1,
\ldots,\theta^*_D$ such that
$A^* = \sum_{i=0}^D \theta^*_iE^*_i$.
Observe
$A^*E^*_i = E^*_iA^* =  \theta^*_iE^*_i$ for $0 \leq i \leq D$.
By \cite[Lemma 3.11]{terwSub1} the 
scalars $\theta^*_0,\theta^*_1,\ldots,\theta^*_D$ are in $\R.$ 
The scalars $\theta^*_0,\theta^*_1,\ldots,\theta^*_D$ are mutually
distinct 
since $A^*$ generates $M^*$. We call $\theta^*_i$ the {\it dual eigenvalue}
of $\Gamma$ associated with $E^*_i$ $(0 \leq i\leq D)$.

\medskip
\noindent 
We recall the subconstituents of $\Gamma $.
From
(\ref{DEFDEI}) we find
\begin{equation}\label{DEIV}
E_i^*V = \mbox{span}\{\hat{y} \;|\; y \in X, \quad \partial(x,y)=i\}
\qquad (0 \le i \le D).
\end{equation}
By 
(\ref{DEIV})  and since
 $\{\hat{y}\;|\;y \in X\}$ is an orthonormal basis for $V$
 we find
\begin{eqnarray*}
\label{vsub}
V = E_0^*V+E_1^*V+ \cdots +E_D^*V \qquad \qquad 
{\rm (orthogonal\ direct\ sum}).
\end{eqnarray*}
For $0 \leq i \leq D$ the space $E^*_iV$ is the eigenspace
of $A^*$ associated with $\theta^*_i$.
We call $E_i^*V$ the {\it $i$th subconstituent} of $\Gamma$
with respect to $x$.

\medskip
\noindent
We recall the subconstituent algebra of $\Gamma $.
Let $T=T(x)$ denote the subalgebra of $\hbox{Mat}_X(\C)$ generated by 
$M$ and $M^*$. 
We call $T$ the {\it subconstituent algebra} 
(or {\it Terwilliger algebra}) of $\Gamma$ 
 with respect to $x$ \cite[Definition 3.3]{terwSub1}.
We observe $T$ is generated by $A$ and $A^*$.
We observe $T$ has finite dimension. Moreover $T$ is 
semi-simple since it
is closed under the conjugate transponse map
\cite[p.~157]{CR}. 
See
\cite{curtin1,
curtin2,
curtin6,
egge1,
go,
go2,
hobart,
tanabe,
terwSub1,
terwSub2,
terwSub3}
for more information on the subconstituent
algebra.

\medskip
\noindent 
For the rest of this paper we adopt the
following notational convention.

\begin{definition}
\label{setup}
\rm
We assume $\Gamma=(X,R)$ is a  distance-regular graph with
diameter $D\geq 3$. We assume $\Gamma$ is $Q$-polynomial
with respect to the ordering $E_0, E_1, \ldots, E_D$
of the primitive idempotents. We fix $x \in X$ 
and write $A^*=A^*(x)$, 
$E^*_i=E^*_i(x)$ $(0 \leq i \leq D)$,
$T=T(x)$. We abbreviate $V=\C^X$.
For notational convenience we define
$E_{-1}=0$, 
$E_{D+1}=0$ and 
$E^*_{-1}=0$, 
$E^*_{D+1}=0$.
\end{definition}

\noindent We have some comments.

\begin{lemma}
\label{lem:incl}
\cite[Lemma 3.2]{terwSub1}
With reference to Definition \ref{setup},
the following (i), (ii) hold.
\begin{enumerate}
\item
$AE^*_iV \subseteq E^*_{i-1}V +E^*_iV+E^*_{i+1}V
$ $(0 \leq i \leq D)$.
\item $A^*E_iV \subseteq E_{i-1}V +E_iV+E_{i+1}V
$ $(0 \leq i \leq D)$.
\end{enumerate}
\end{lemma}

\begin{lemma}
\label{lem:incl2}
With reference to Definition \ref{setup},
the following (i)--(iv) hold.
\begin{enumerate}
\item $A \sum_{h=0}^iE^*_hV \subseteq \sum_{h=0}^{i+1} E^*_hV
$ $(0 \leq i \leq D)$.
\item $(A-\theta_i I)\sum_{h=0}^iE_hV =\sum_{h=0}^{i-1} E_hV
$ $(0 \leq i \leq D)$.
\item $A^* \sum_{h=0}^iE_hV \subseteq \sum_{h=0}^{i+1} E_hV
$ $(0 \leq i \leq D)$.
\item $(A^*-\theta^*_i I)\sum_{h=0}^iE^*_hV =\sum_{h=0}^{i-1} E^*_hV
$ $(0 \leq i \leq D)$.
\end{enumerate}
\end{lemma}
{\it Proof:}
(i) Immediate from Lemma
\ref{lem:incl}(i).
\\
\noindent (ii) Recall $AE_j=\theta_j E_j$ for $0 \leq j \leq D$.
\\
\noindent (iii) Immediate from Lemma
\ref{lem:incl}(ii).
\\
\noindent (iv) Recall $A^*E^*_j=\theta^*_j E^*_j$
for $0 \leq j \leq D$.
\hfill $\Box$ \\

\section{The irreducible $T$-modules}

\noindent In this section 
we recall some results
on $T$-modules for later use.

\medskip
\noindent
With reference to Definition \ref{setup},
by a {\it T-module}
we mean a subspace $W \subseteq V$ such that $BW \subseteq W$
for all $B \in T.$ 
 Let $W$ denote a $T$-module. Then $W$ is said
to be {\it irreducible} whenever $W$ is nonzero and $W$ contains 
no $T$-modules other than 0 and $W.$ Let $W, W^\prime$ denote
$T$-modules. By an {\it isomorphism of $T$-modules}
from $W$ to $W^\prime$ we
mean an isomorphism of vector spaces
$\sigma: W \rightarrow W^\prime$
such that
$(\sigma B- B \sigma)W = 0$ for
all $B \in T$.
The modules $W, W^\prime$ are said to be {\it isomorphic as $T$-modules}
whenever
there exists an isomorphism of $T$-modules from $W$ to $W^\prime.$
Any two nonisomorphic
irreducible $T$-modules are orthogonal \cite[Lemma 3.3]{curtin1}.

\medskip
\noindent
Let $W$ denote a $T$-module and let 
$W'$ denote a  
$T$-module contained in $W$.
Then the orthogonal complement of $W'$ in $W$ is a $T$-module 
\cite[p.~802]{go2}.
It follows that each $T$-module
is an orthogonal direct sum of irreducible $T$-modules.
In particular $V$ is an orthogonal direct sum of irreducible $T$-modules.

\medskip
\noindent 
Let $W$ denote an irreducible $T$-module.
By the {\it endpoint} of $W$ we mean
$\mbox{min}\lbrace i |0\leq i \leq D, \; E^*_iW\not=0\rbrace $.
By the {\it diameter} of $W$ we mean
$ |\lbrace i | 0 \leq i \leq D,\; E^*_iW\not=0 \rbrace |-1 $.
By the {\it dual endpoint} of $W$ we mean
$\mbox{min}\lbrace i |0\leq i \leq D, \; E_iW\not=0\rbrace $.
By
the {\it dual diameter} of $W$ we mean
$ |\lbrace i | 0 \leq i \leq D,\; E_iW\not=0 \rbrace |-1 $.
The diameter of $W$ is  equal to the dual diameter of
$W$
\cite[Corollary 3.3]{aap1}.
There exists a unique irreducible $T$-module with diameter $D$.
We call this module the {\it primary} $T$-module.
The primary $T$-module 
has basis $A_0{\hat x}, \ldots, A_D{\hat x}$
\cite[Lemma 3.6]{terwSub1}.

\begin{lemma}
\cite[Lemma 3.4, Lemma 3.9, Lemma 3.12]{terwSub1}
\label{lem:basic}
With reference to Definition \ref{setup},
let $W$ denote an irreducible $T$-module with endpoint $\rho$,
dual endpoint $\tau$, and diameter $d$.
Then $\rho,\tau,d$ are nonnegative integers such that $\rho+d\leq D$ and
$\tau+d\leq D$. Moreover the following (i)--(iv) hold.
\begin{enumerate}
\item 
$E^*_iW \not=0$ if and only if $\rho \leq i \leq \rho+d$, 
$ \quad (0 \leq i \leq D)$.
\item
$W = \sum_{h=0}^{d} E^*_{\rho+h}W \qquad (\mbox{orthogonal direct sum}). $
\item 
$E_iW \not=0$ if and only if $\tau \leq i \leq \tau+d$,
$ \quad (0 \leq i \leq D)$.
\item
$W = \sum_{h=0}^{d} E_{\tau+h}W \qquad (\mbox{orthogonal direct sum}). $
\end{enumerate}
\end{lemma}

\begin{lemma}
\label{lem:wincl}
With reference to Definition \ref{setup},
let $W$ denote an irreducible $T$-module with
endpoint $\rho$,
dual endpoint $\tau$, and diameter $d$.
Then the following (i), (ii) hold.
\begin{enumerate}
\item
$AE^*_{\rho+i}W \subseteq E^*_{\rho+i-1}W +E^*_{\rho+i}W+E^*_{\rho+i+1}W
$ $ \qquad (0 \leq i \leq d)$.
\item $A^*E_{\tau+i}W \subseteq E_{\tau+i-1}W +E_{\tau+i}W+E_{\tau+i+1}W
$ $ \qquad (0 \leq i \leq d)$.
\end{enumerate}
\end{lemma}
\noindent 
{\it Proof:} (i) Follows from 
Lemma
\ref{lem:incl}(i) and since $E^*_jW=E^*_jV \cap W$ for $0 \leq j \leq D$.
\\
\noindent (ii) Follows from 
Lemma
\ref{lem:incl}(ii) and since $E_jW=E_jV \cap W$ for $0 \leq j \leq D$.
\hfill $\Box$ \\

\begin{remark} 
\label{rem:tdpair}
\rm
With reference to Definition \ref{setup},
let $W$ denote an irreducible $T$-module.
Then $A$ and  $A^*$ act on $W$ as a tridiagonal pair
in the sense of \cite[Definition 1.1]{itt}.
This follows from
Lemma \ref{lem:basic},
Lemma 
\ref{lem:wincl}, and since $A,A^*$ together generate
$T$. See 
\cite{shape,
LS99,
qSerre,
LS24,
conform,
lsint,
Terint,
TLT:split,
TLT:array,
qrac,
awrel} for information on tridiagonal pairs. 
\end{remark}

\begin{lemma}
\cite[Lemma 5.1, Lemma 7.1]{caugh2}
\label{lem:basineq}
With reference to Definition \ref{setup},
let $W$ denote an irreducible $T$-module with endpoint
$\rho$, dual endpoint $\tau$, and diameter $d$. Then 
the following (i), (ii) hold.
\begin{enumerate}
\item $2\rho+d\geq D$.
\item $2\tau+d\geq D$.
\end{enumerate}
\end{lemma}

\begin{lemma}
\label{thm:wsplit}
With reference to Definition \ref{setup}, let $W$ 
denote an irreducible $T$-module with endpoint $\rho$,
dual endpoint $\tau$, and diameter $d$.
Then 
\begin{eqnarray}
\label{eq:splitd}
W = \sum_{h=0}^d W_h \qquad (\mbox{direct sum}),
\end{eqnarray}
where
\begin{eqnarray}
W_h = (E^*_\rho W+\cdots+ E^*_{\rho+h}W)\cap (E_\tau W+\cdots + E_{\tau+d-h}W)
\qquad (0 \leq h\leq d).
\label{eq:whdef}
\end{eqnarray}
\end{lemma}
\noindent 
{\it Proof:} Immediate from 
Remark
\ref{rem:tdpair} and
\cite[Theorem 4.6]{itt}.
\hfill $\Box$ \\

\begin{remark}
\rm
The sum
(\ref{eq:splitd}) is not orthogonal in general.
\end{remark}

\section{The displacement decomposition}

\noindent In this section we introduce the
{\it displacement decomposition} for the standard module.

\begin{definition}
\label{def:shift}
\rm
With reference to Definition
 \ref{setup}, let $W$ denote an irreducible $T$-module.
By the {\it displacement} of $W$ we mean the integer $\rho+\tau+d-D$, where 
$\rho,\tau,d$ denote respectively 
the endpoint, dual endpoint, and diameter of $W$.
\end{definition}

\begin{lemma}
\label{lem:shiftpos}
With reference to Definition
\ref{setup}, let $W$ denote an irreducible $T$-module
with displacement $\eta$. Then $0 \leq \eta \leq D$.
\end{lemma}
\noindent 
{\it Proof:} Let $\rho,\tau,d$ denote respectively 
the endpoint, dual endpoint,
and diameter of $W$.
By Lemma
\ref{lem:basineq}
we have $2\rho+d\geq D$ and $2\tau+d\geq D$; adding these inequalities
we find $\rho+\tau+d\geq D$ so 
$\eta\geq 0$. By 
Lemma \ref{lem:basic}
we have $\rho\leq D$ and $\tau+d\leq D$. Combining these
inequalities
we find $\rho+\tau+d\leq 2D$ 
so
$\eta\leq D$.
\hfill $\Box$ \\

\begin{definition}
\label{def:modl}
\rm
With reference to Definition \ref{setup},
For $0 \leq \eta \leq D$ we let  ${V}_\eta$ denote
the subspace of $V$ spanned by the irreducible $T$-modules
that have displacement $\eta$.
We observe $V_\eta$ is a $T$-module.
\end{definition}

\begin{lemma}
\label{lem:vsorthog}
With reference to Definition
\ref{setup},
\begin{eqnarray}
V = \sum_{\eta=0}^D {V}_\eta \qquad (\mbox{orthogonal direct sum}).
\label{eq:disp}
\end{eqnarray}
\end{lemma}
\noindent {\it Proof:}
We mentioned earlier that
$V$ is spanned by the irreducible $T$-modules.
By  
Lemma
\ref{lem:shiftpos}
and 
Definition
\ref{def:modl},
 each of these modules
is contained in one of $V_0, V_1, \ldots, V_D$.
Therefore
$V=\sum_{\eta=0}^D V_\eta$.
To show this sum is orthogonal and direct,
it suffices to show
$V_0, V_1, \ldots, V_D$ are mutually orthogonal.
For distinct integers $i,j$ $(0 \leq i,j \leq D)$
observe $V_i, V_j$ are orthogonal since
the isomorphism classes of irreducible
$T$-modules that span $V_i$ 
are distinct from the isomorphism classes of irreducible $T$-modules
that span 
$V_j$.
We have now shown
$V_0, V_1, \ldots, V_D$ are mutually orthogonal so
the sum
$\sum_{\eta=0}^D V_\eta$ is orthogonal and direct.
\hfill $\Box$ \\

\begin{definition}
\rm
We call the sum
(\ref{eq:disp}) the 
{\it displacement decomposition} of $V$ with respect to $x$.
\end{definition}

\section{The split decomposition}

\noindent In this section we introduce the {\it split decomposition}
of the standard module.

\begin{definition}
\label{uij}
With reference to Definition \ref{setup},
for $-1 \leq i,j\leq D$ we define
\begin{eqnarray}
\label{eq:vij}
V_{ij} = (E^*_0V+E^*_1V+\cdots+E^*_iV)\cap (E_0V+E_1V+\cdots +E_jV).
\end{eqnarray}
We observe $V_{ij}=0$ if $i=-1$ or $j=-1$.
\end{definition}

\noindent In the following three lemmas
we make some observations concerning Definition
\ref{uij}. In each case the proof is 
routine and omitted.

\begin{lemma}
\label{uobs1}
With reference to Definition \ref{setup},
for $0\leq i,j\leq D$ the space $V_{ij}$
consists of those vectors $v \in V$ such
that 
 $E^*_h v=0$ for $i< h \leq D$ and 
 $E_h v=0$ for $j<h \leq D$.
\end{lemma}

\begin{lemma}
\label{lem:uobs2}
With reference to Definition
\ref{setup},  we have
$V_{i-1,j}\subseteq V_{ij}$ 
and $V_{i,j-1}\subseteq V_{ij}$ for $0 \leq i,j\leq D$.
\end{lemma}

\begin{lemma}
\label{lem:uobs4}
With reference to Definition
\ref{setup}, the following (i)--(iii) hold.
\begin{enumerate}
\item $V_{iD}  = E^*_0V+E^*_1V+\cdots +E^*_iV $ $(0 \leq i \leq D)$.
\item $V_{Dj}  = E_0V+E_1V+\cdots +E_jV $ $(0 \leq j \leq D)$.
\item $V_{DD}=V$.
\end{enumerate}
\end{lemma}

\noindent Later in the paper 
we will show
$V_{ij}=0$ 
if $i+j<D$, $(0 \leq i,j\leq D)$.

\begin{definition}
\rm
\label{def:tildeu}
With reference to Definition
\ref{setup},
for $0 \leq i,j\leq D$ we let  ${\tilde V}_{ij}$ denote
the orthogonal complement of $V_{i,j-1}+V_{i-1,j}$ in
$V_{ij}$. For notational convenience we define
${\tilde V}_{ij}:=0$ unless $i,j \in \lbrace 0,1,\ldots, d\rbrace $.
\end{definition}

\noindent Our next goal is to show
$V_{rs}=\sum_{i=0}^r \sum_{j=0}^s 
{\tilde V}_{ij}$ (direct sum) for
$0 \leq r,s\leq D$. 
We will use the following lemma.

\begin{lemma}
\label{lem:dimsum}
With reference to Definition
\ref{setup},
\begin{eqnarray}
\label{eq:dimsum}
\mbox{dim}\,{\tilde V}_{ij}  =  
\mbox{dim}\,V_{ij} - 
\mbox{dim}\,V_{i,j-1} -  
\mbox{dim}\,V_{i-1,j} + 
\mbox{dim}\,V_{i-1,j-1}
\end{eqnarray}
 for $0 \leq i,j\leq D$. 
\end{lemma}
\noindent {\it Proof:}
Let $z$ denote the
dimension of
$V_{i,j-1}+V_{i-1,j}$.
The space ${\tilde V}_{ij}$ is the orthogonal complement of
$V_{i,j-1}+V_{i-1,j}$ in 
$V_{ij}$ so
$\mbox{dim}\,{\tilde V}_{ij} + z = \mbox{dim}\, V_{ij}$.
Using Definition
\ref{uij}
we find
$V_{i,j-1}\cap V_{i-1,j} = V_{i-1,j-1}$
so  
$z+\mbox{dim}\,V_{i-1,j-1}=
\mbox{dim}\,V_{i,j-1}+\mbox{dim}\,V_{i-1,j}$.
From these comments we routinely obtain 
(\ref{eq:dimsum}).
\hfill $\Box$ \\

\begin{theorem}
\label{thm:usplitdec}
With reference to Definition
\ref{setup}, for $0 \leq r,s\leq D$ we have
\begin{eqnarray*}
V_{rs} =\sum_{i=0}^r
\sum_{j=0}^s
{\tilde V}_{ij}
\qquad (\mbox{direct sum}).
\end{eqnarray*}
\end{theorem}
\noindent {\it Proof:}
\noindent  
We first show
\begin{eqnarray}
V_{rs} =\sum_{i=0}^r
\sum_{j=0}^s
{\tilde V}_{ij}.
\label{eq:justsum}
\end{eqnarray}
The proof is by induction on $r+s$. The result is trivial for
$r+s=0$ so assume $r+s>0$.
Recall 
${\tilde V}_{rs}$ is the orthogonal complement of
$V_{r,s-1} 
+V_{r-1,s}
$
in $V_{rs}$. Therefore
\begin{eqnarray}
\label{eq:sum1}
V_{rs} ={\tilde V}_{rs}
+
V_{r,s-1} 
+V_{r-1,s}.
\end{eqnarray}
By induction we have both
\begin{eqnarray}
\label{eq:sum2}
V_{r,s-1} 
=
\sum_{i=0}^r
\sum_{j= 0}^{s-1}
{\tilde V}_{ij},
\qquad \qquad 
V_{r-1,s} 
=
\sum_{i=0}^{r-1}
\sum_{j=0}^s
{\tilde V}_{ij}.
\end{eqnarray}
Combining 
(\ref{eq:sum1}), (\ref{eq:sum2})
 we routinely obtain  (\ref{eq:justsum}).
We now show the sum
  (\ref{eq:justsum}) is direct.
From Lemma
\ref{lem:dimsum} we routinely obtain
\begin{eqnarray*}
\mbox{dim}\,V_{rs} =\sum_{i=0}^r
\sum_{j=0}^s
\mbox{dim}\,{\tilde V}_{ij}
\end{eqnarray*}
and it follows 
  the sum (\ref{eq:justsum}) is direct.
\hfill $\Box$ \\

\begin{corollary}
\label{cor:uspdec}
With reference to Definition \ref{setup},
\begin{eqnarray}
V =\sum_{i=0}^D
\sum_{j=0}^D
{\tilde V}_{ij}
\qquad (\mbox{direct sum}).
\label{eq:vsplt}
\end{eqnarray}
\end{corollary}
\noindent {\it Proof:}
Set $r=D$ and $s=D$ in Theorem
\ref{thm:usplitdec} and use
Lemma
\ref{lem:uobs4}(iii).
\hfill $\Box$ \\

\begin{definition}
\label{def:splitdv}
\rm
We call the sum 
(\ref{eq:vsplt})
the {\it split decomposition} of $V$ with respect to
 $x$. This decomposition is
not orthogonal in general.
\end{definition}

\section{The displacement and split decompositions}

\noindent In this section we describe the relationship
between the displacement
decomposition and the 
split decomposition.
Our main result is the following.
With reference to Definition \ref{setup},
for $0 \leq \eta \leq D$ we show $V_\eta = 
\sum {\tilde V}_{ij}$,
where the sum is over all ordered pairs $i,j$
such that $0 \leq i,j\leq D$ and $i+j=D+\eta$.
We begin with a lemma.

\begin{lemma}
\label{lem:pre}
With reference to Definition \ref{setup},
let $W$ denote an irreducible $T$-module with endpoint $\rho$,
dual endpoint $\tau$, and diameter $d$.
Let the subspaces $W_0, W_1, \ldots, W_d$
be as in Lemma
\ref{thm:wsplit}. Then 
$W_h\subseteq {\tilde V}_{\rho+h,\tau+d-h}$ for $0 \leq h \leq d$.
\end{lemma}
\noindent {\it Proof:}
Comparing 
(\ref{eq:whdef}) and
(\ref{eq:vij}) 
we  find 
$W_h \subseteq V_{\rho+h,\tau+d-h}$.
We  show $W_h$ is orthogonal to 
$V_{\rho+h-1,\tau+d-h}+
V_{\rho+h,\tau+d-h-1}$.
For $w \in W_h$ and for 
$v \in 
V_{\rho+h-1,\tau+d-h}$ we show $\langle w,v\rangle =0$.
Let $W^\perp$ denote the orthogonal complement of $W$ in
$V$. Observe $V=W + W^\perp$ (direct sum) and that
$W^\perp$ is a $T$-module.
Observe there exists $w_1 \in W$ and $v_1 \in W^\perp$ such
that $v=w_1+v_1$.
By the construction $w \in W$ and $v_1 \in W^\perp$ so
$\langle w,v_1\rangle =0$.
We show $w_1=0$.
By Lemma
\ref{uobs1} and
since $v \in 
V_{\rho+h-1,\tau+d-h}$ 
we find $E^*_iv=0$ for $\rho+h\leq i \leq D$ and
 $E_jv=0$ for $\tau+d-h+1\leq j \leq D$.
Since
$V=W + W^\perp$ is a direct sum of $T$-modules
we find $E^*_iw_1=0$ for $\rho+h\leq i \leq D$ and
 $E_jw_1=0$ for $\tau+d-h+1\leq j \leq D$.
Since $w_1 \in W$ and since $W$ has endpoint $\rho$
we have $E^*_iw_1=0$ for $0 \leq i \leq \rho -1$.
Similarly since $W$ has dual endpoint $\tau$ we have
$E_jw_1=0$ for $0 \leq j \leq \tau-1$.
From these comments we find
\begin{eqnarray}
w_1 \in (E^*_\rho W+\cdots + E^*_{\rho+h-1} W)\cap
(E_\tau W+\cdots +E_{\tau + d-h} W).
\label{eq:preint}
\end{eqnarray}
Using (\ref{eq:whdef}) we find the intersection on the
right in
(\ref{eq:preint}) is equal to $W_h\cap W_{h-1}$, where $W_{-1}=0$.
The sum
(\ref{eq:splitd}) is direct so
$W_h \cap W_{h-1}=0$. We now see $w_1=0$.
Now $v=v_1$ so $\langle w, v\rangle =0$.
We have now shown 
 $W_h$ is orthogonal to 
$V_{\rho+h-1,\tau+d-h}$. By a similar argument
we find 
 $W_h$ is orthogonal to 
$V_{\rho+h,\tau+d-h-1}$.
We conclude $W \subseteq
{\tilde V}_{\rho+h,\tau+d-h}$.
\hfill $\Box$ \\

\begin{theorem}
\label{thm:connect}
With reference to Definition \ref{setup},
the following (i)--(iii) hold.
\begin{enumerate}
\item
For $0 \leq \eta \leq D$ we have
$V_\eta = \sum {\tilde V}_{ij}$, where the sum is over all ordered
pairs $i,j$ such that $0 \leq i,j\leq D$ and $i+j=D+\eta$.
\item ${\tilde V}_{ij}=0 $ if $i+j<D$, $\qquad (0 \leq i,j\leq D)$.
\item $V_{ij}=0 $ if $i+j<D$, $ \qquad (0 \leq i,j\leq D)$.
\end{enumerate}
\end{theorem}
\noindent {\it Proof:}
(i), (ii) 
For $-D \leq \eta \leq D$ we define
$V'_\eta = \sum {\tilde V}_{ij}$ where the sum is over all ordered
pairs $i,j$ such that $0 \leq i,j\leq D$ and $i+j=D+\eta$.
Using 
(\ref{eq:vsplt}) we find
$V=\sum_{\eta = -D}^D V'_\eta $ (direct sum).
We show $V'_\eta = 0 $ for $-D \leq \eta < 0$ 
and 
$V'_\eta=V_\eta $ for $0 \leq \eta \leq D$.
Since the sums 
$V=\sum_{\eta = 0}^D V_\eta $
and
$V=\sum_{\eta = -D}^D V'_\eta $ 
 are direct it suffices 
 to show
$V_\eta \subseteq V'_\eta$ for $0 \leq \eta \leq D$.
Let $\eta $ be given.
Let $W$ denote an irreducible $T$-module with
displacement $\eta$. 
Combining Lemma
\ref{thm:wsplit}
and Lemma
\ref{lem:pre}
we find $W\subseteq V'_\eta$.
The space $V_\eta$ is spanned by the irreducible $T$-modules
that have displacement $\eta$; therefore
$V_\eta \subseteq V'_\eta$.
We have now shown
$V_\eta \subseteq V'_\eta$ for $0 \leq \eta \leq D$.
We conclude $V'_\eta = 0 $ for $-D \leq \eta < 0$ 
and 
$V'_\eta=V_\eta $ for $0 \leq \eta \leq D$.
Lines (i), (ii) follow.
\\
\noindent (iii) Combine (ii) above with
Theorem
\ref{thm:usplitdec}.
\hfill $\Box$ \\

\noindent We have some comments.

\begin{theorem}
\label{thm:interpre}
With reference to Definition
\ref{setup},  
for $0 \leq i,j\leq D$ such that $i+j\geq D$,
and for $0 \leq  \eta \leq D$, 
\begin{eqnarray*}
V_{ij}\cap V_\eta=
\sum {\tilde V}_{rs},
\end{eqnarray*}
where the sum is over all ordered pairs
$r,s$ such that $0 \leq r\leq i$  and
$0 \leq s\leq j$ and $r+s-D=\eta$.
\end{theorem}
\noindent {\it Proof:}
Combine
Theorem \ref{thm:usplitdec}
and Theorem
\ref{thm:connect}(i).
\hfill $\Box$ \\

\begin{corollary}
\label{thm:inter}
With reference to Definition
\ref{setup},  
for $0 \leq i,j\leq D$ such that $i+j\geq D$,
we have ${\tilde V}_{ij} = V_{ij}\cap V_\eta$ where
$\eta = i+j-D$.
\end{corollary}
\noindent {\it Proof:}
Apply Theorem
\ref{thm:interpre} with $\eta = i+j-D$.
\hfill $\Box$ \\

\section{The action of $A$ and $A^*$ on the split decomposition}

In this section we describe
how the adjacency matrix and the dual adjacency matrix
 act on the split
decomposition.

\begin{theorem}
With reference to Definition
\ref{setup}, the following (i), (ii) hold.
\begin{enumerate}
\item $(A-\theta_jI){\tilde V}_{ij} \subseteq {\tilde V}_{i+1,j-1}
$ 
$(0 \leq i,j\leq D)$.
\item $(A^*-\theta^*_iI){\tilde V}_{ij} \subseteq {\tilde V}_{i-1,j+1}
$ 
$(0 \leq i,j\leq D)$.
\end{enumerate}
\end{theorem}
\noindent {\it Proof:}
(i) Assume $i+j\geq D$; otherwise ${\tilde V}_{ij}=0$ 
and the result is trivial.
For convenience we treat the cases $i=D$ and $i<D$ separately.
To obtain the result for the case $i=D$,
we show $(A-\theta_jI){\tilde V}_{Dj}=0$.
From Corollary 
\ref{thm:inter} (with $i=D$ and $\eta =j$)  we have
${\tilde V}_{Dj}=V_{Dj}\cap V_j$.
Using
Lemma
\ref{lem:incl2}(ii) and 
Lemma
\ref{lem:uobs4}(ii)
we find
 $(A-\theta_jI)V_{Dj}=V_{D,j-1}$.
Therefore
 $(A-\theta_jI){\tilde V}_{Dj}\subseteq V_{D,j-1}$.
Recall $V_j$ is a $T$-module so
 $(A-\theta_jI)V_j\subseteq V_j$.
Therefore
 $(A-\theta_jI){\tilde V}_{Dj}\subseteq V_j$.
Now
\begin{eqnarray*}
 (A-\theta_jI){\tilde V}_{Dj}&\subseteq&V_{D,j-1}\cap V_j
\\
&=& 0
\end{eqnarray*}
in view of 
 Theorem
\ref{thm:interpre}.
We have now shown $(A-\theta_jI){\tilde V}_{Dj}=0$
so we are done for the case $i=D$.
Next assume $i<D$.
From Corollary 
\ref{thm:inter} we have
${\tilde V}_{ij}=V_{ij}\cap V_\eta$ where $\eta = i+j-D$.
Using 
Lemma 
\ref{lem:incl2}
and 
(\ref{eq:vij}) we find
$(A-\theta_jI)V_{ij}\subseteq V_{i+1,j-1}$.
Therefore 
$(A-\theta_jI){\tilde V}_{ij}\subseteq V_{i+1,j-1}$.
Recall $V_\eta$ is a $T$-module so
$(A-\theta_jI)V_\eta \subseteq V_\eta $.
Therefore 
$(A-\theta_jI){\tilde V}_{ij} \subseteq V_\eta $.
Now
\begin{eqnarray*}
(A-\theta_jI){\tilde V}_{ij} &\subseteq& 
 V_{i+1,j-1}
\cap V_\eta 
\\
&=& {\tilde V}_{i+1,j-1}
\end{eqnarray*}
in view of Corollary 
\ref{thm:inter}.
\\
\noindent (ii) Similar to the proof of (i) above.
\hfill $\Box$ \\

\section{An application}

\noindent In this section we give an application of
Theorem
\ref{thm:connect}(iii).
We first give two definitions.

\begin{definition}
\rm
Let $\Gamma=(X,R)$ denote a distance-regular graph 
 with standard module $V$.
For $v \in V$, by the {\it support} of $v$
we mean the subset of $X$ consisting
of those vertices $y$ such that coordinate $y$ of $v$ is nonzero.
\end{definition}

\begin{definition}
\label{def:dualwidth}
\rm
\cite[Section 4]{wdw}
Let $\Gamma$ denote a distance-regular graph with diameter
$D\geq 3$. Assume $\Gamma$ is $Q$-polynomial
with respect to the ordering $E_0, E_1, \ldots, E_D$ of 
the primitive idempotents.
Let $v$ denote a nonzero vector in the standard module $V$. 
By the {\it dual width} of $v$ we mean
\begin{eqnarray*}
  \mbox{max}\lbrace i |0 \leq i \leq D, \;E_iv\not=0\rbrace.
\end{eqnarray*}
\end{definition}

\begin{theorem}
\label{thm:mth}
Let $\Gamma=(X,R)$ denote a distance-regular graph with diameter
$D\geq 3$. Assume $\Gamma$ is $Q$-polynomial
with respect to the ordering $E_0, E_1, \ldots, E_D$ of 
the primitive idempotents.
Let $v$ denote a nonzero vector in the standard module $V$ and let $g$ denote
the corresponding dual width from Definition
\ref{def:dualwidth}.
Then for  all $x \in X$ there exists
$y$ in the support of $v$ such that
\begin{eqnarray}
\partial (x,y)\geq D-g.
\label{eq:mineq}
\end{eqnarray}
\end{theorem}
{\it Proof:}
We assume the result is false and obtain  a contradiction.
By this assumption there exists $x \in X$ such that
$\partial(x,y)<D-g$ for all vertices $y$ in the support of $v$.
Abbreviate $E^*_i=E^*_i(x)$ for $0 \leq i \leq D$.
Then $v \in E^*_0V+\cdots + E^*_fV$ where $f=D-g-1$.
Using Definition
\ref{def:dualwidth} we find
 $v \in E_0V+\cdots + E_{g}V$.
Now
\begin{eqnarray*}
v &\in&
(E^*_0V+\cdots + E^*_fV)\cap
(E_0V+\cdots + E_{g}V)
\\
&=&V_{fg}.
\end{eqnarray*}
We mentioned $f=D-g-1$ so  $f+g<D$;
combining this with
Theorem
\ref{thm:connect}(iii) 
we find
$V_{fg}=0$. Now $v=0$ for a contradiction.
The result follows.
\hfill $\Box$ \\

\begin{remark}
\rm
Referring to Theorem
\ref{thm:mth}, pick any $x \in X$. 
If $v$ is not orthogonal to the primary
module for $T(x)$ then
(\ref{eq:mineq}) follows
from
\cite[Equation (2.8)]{roos}. See also \cite[Lemma 1]{wdw}. 
\end{remark}

\section{Directions for further research}

\noindent In this section we give some suggestions for
further research.

\begin{problem} \rm
With reference to
Definition
\ref{setup},
recall that for $0 \leq i,j\leq D$ the space ${\tilde V}_{ij}$
depends on $x$.
Does the dimension of 
 ${\tilde V}_{ij}$
depend on $x$?
\end{problem}

\begin{problem} \rm
With reference to
Definition
\ref{setup},
let $W$ denote an irreducible $T$-module and
consider the multiplicity with which $W$ appears in
$V$.
In general this multiplicity is not determined by 
the intersection numbers of $\Gamma$ 
\cite{tanabe}.
Is this multiplicity determined by
the intersection numbers of $\Gamma$ and the scalars
$\lbrace \mbox{dim}{\tilde V}_{ij} \;|\;0 \leq i,j\leq D\rbrace $?
\end{problem}

\begin{problem} \rm Let $\Gamma$ denote a 
$Q$-polynomial distance-regular graph. In many cases
$\Gamma$
exists on the top fiber of a ranked poset
\cite{delsarte},
\cite{terqmat},
\cite{terunif}.
For this case
investigate the relationship
between
 the
poset structure and
 the 
split decomposition of $\Gamma$.
\end{problem}

\bigskip
\noindent{\Large\bf Acknowledgements}

\medskip
\noindent 
The author would like to thank 
Brian Curtin, Eric Egge, Michael Lang, Stefko Miklavic,
and Arlene Pascasio 
for giving this manuscript a careful reading and 
offering many valuable suggestions.

\noindent Paul Terwilliger \\
 Department of Mathematics \\ 
University of Wisconsin \\
Van Vleck Hall \\
480 Lincoln Drive \\
Madison Wisconsin \\
USA  53706-1388 \\
Email: terwilli@math.wisc.edu

\end{document}